\documentclass[twoside]{article}
\usepackage{amssymb}
\usepackage{amsmath,amssymb,amsthm}

\begin{document}

\frenchspacing

\pagestyle{myheadings}

\markboth{Dimitar Mekerov}{Connection with skew-symmetric torsion
on almost Norden manifolds}

\numberwithin{equation}{section}
\newtheorem{thm}{Theorem}[section]
\newtheorem{lem}[thm]{Lemma}
\newtheorem{prop}[thm]{Proposition}
\newtheorem{cor}[thm]{Corollary}
\newtheorem{probl}[thm]{Problem}

\newtheorem{defn}{Definition}[section]
\newtheorem{rem}{Remark}[section]
\newtheorem{exa}{Example}



\newcommand{\X}{\mathfrak{X}}
\newcommand{\B}{\mathcal{B}}
\newcommand{\s}{\mathfrak{S}}
\newcommand{\g}{\mathfrak{g}}
\newcommand{\W}{\mathcal{W}}
\newcommand{\Lgr}{\mathrm{L}}
\newcommand{\dd}{\mathrm{d}}

\newcommand{\diag}{\mathrm{diag}}
\newcommand{\End}{\mathrm{End}}
\newcommand{\im}{\mathrm{Im}}
\newcommand{\Id}{\mathrm{Id}}

\newfont{\w}{msbm9 scaled\magstep1}
\def\R{\mbox{\w R}}
\newcommand{\norm}[1]{\left\Vert#1\right\Vert ^2}
\newcommand{\nN}{\norm{N}}
\newcommand{\ad}{{\rm ad}}
\newcommand{\n}{\nabla}
\newcommand{\nJ}{\norm{\nabla J}}

\newcommand{\be}[1]{\begin{equation}\label{#1}}
\newcommand{\ee}{\end{equation}}

\newcommand{\thmref}[1]{Theorem~\eqref{#1}}
\newcommand{\propref}[1]{Proposition~\eqref{#1}}
\newcommand{\secref}[1]{\S\eqref{#1}}
\newcommand{\lemref}[1]{Lemma~\eqref{#1}}
\newcommand{\dfnref}[1]{Definition~\eqref{#1}}

\frenchspacing

\hyphenation{Her-mi-ti-an ma-ni-fold} \hyphenation{iso-tro-pic
mani-fold mani-folds}


\title{On the geometry of the connection with totally skew-symmetric torsion 
on almost complex manifolds 
with Norden metric}

\author{Dimitar Mekerov}

\maketitle

\noindent University of Plovdiv, Faculty of Mathematics and
Informatics, 236 Bulgaria Blvd., 4003 Plovdiv, Bulgaria,
mircho@uni-plovdiv.bg

\bigskip

\noindent
\textbf{Abstract} \\
We consider an almost complex manifold with Norden metric (i.~e. a
metric with respect to which the almost complex structure is an
anti-isometry). On such a manifold we study a linear connection
preserving the almost complex structure and the metric and having
a totally skew symmetric torsion tensor (i.~e. a 3-form). We prove
that if a non-K\"ahler almost complex manifold with Norden metric
admits such connection then the manifold is quasi-K\"ahlerian
(i.~e. has non-integrable almost complex structure). We prove that
this connection is unique, determine its form, and construct an
example of it on a Lie group. We consider the case when the
manifold admits a connection with parallel totally skew-symmetric
torsion and the case when such connection has a K\"ahler curvature
tensor. We get necessary and sufficient conditions for an
isotropic K\"ahler manifold with Norden metric.




\section{Introduction}

There is a strong interest in the metric connections with totally
skew-symmet\-ric torsion tensor (3-form). These connections arise
in a natural way in theoretical and mathematical physics. For
example, such a connection is of particular interest in string
theory \cite{Stro}. In mathematics this connection was used by
Bismut to prove the local index theorem for non-K\"ahler Hermitian
manifolds \cite{Bis}. A connection with totally skew-symmetric
torsion tensor is called a KT connection by physicists,  and among
mathematicians this connection is known as a Bismut connection.

According to Gauduchon \cite{Gaud}, on any non-K\"ahler Hermitian
manifold, there exists a unique Hermitian connection (i.~e. one
preserving the almost complex structure and the metric), whose
torsion tensor is totally skew-symmetric.

The goal of the present work is to solve the analogous problem for
almost complex manifolds with Norden metric. We build upon the
results in \cite{Mek1}, where we have considered a connection with
totally skew-symmetric torsion tensor in the case of
quasi-K\"ahler manifolds with Norden metric.

The main results in the present paper are contained in
\thmref{thm-2.1}. We prove that if a non-K\"ahler almost complex
manifold with Norden metric admits a linear connection $\n'$
preserving the almost complex structure and the metric and having
totally skew-symmetric tensor, then the manifold is
quasi-K\"ahlerian and the connection $\n'$ is unique. We get the
form of $\n'$.

We construct a 4-parametric family of 4-dimensional quasi-K\"ahler
manifolds with Norden metric by a Lie group, and derive the
connection $\n'$ for an arbitrary manifold of this family.

We consider the case when $\n'$ has a parallel torsion and obtain
a relation between the scalar curvatures of this connection and
the Levi-Civita connection. We show that the manifold is isotropic
K\"ahlerian iff these curvatures are equal.

We prove that if the manifold admits a parallel connection with a
K\"ahler curvature tensor, then the manifold is isotropic
K\"ahlerian.


\section{Preliminaries}

Let $(M,J,g)$ be a $2n$-dimensional \emph{almost complex manifold
with Norden metric}, i.~e.
\begin{equation*}
J^2x=-x, \qquad g(Jx,Jy)=-g(x,y),
\end{equation*}
for all differentiable vector fields $x$, $y$ on $M$. The
\emph{associated metric} $\tilde{g}$ of $g$ on $M$, given by
$\tilde{g}(x,y)=g(x,Jy)$, is a Norden metric, too. The signature
of both metrics is necessarily $(n,n)$.

Further, $x$, $y$, $z$, $w$ will stand for arbitrary
differentiable vector fields on $M$ (or vectors in the tangent
space  of $M$ at an arbitrary point $p\in M$).

The Levi-Civita connection of $g$ is denoted by $\nabla$. The
tensor filed $F$ of type $(0,3)$ on $M$ is defined by
\begin{equation}\label{1.1}
F(x,y,z)=g\bigl( \left( \nabla_x J \right)y,z\bigr).
\end{equation}
It has the following properties \cite{GrMeDj}:
\begin{equation*}
F(x,y,z)=F(x,z,y)=F(x,Jy,Jz),\quad F(x,Jy,z)=-F(x,y,Jz).
\end{equation*}

In \cite{GaBo}, the considered manifolds are classified into eight
classes with respect to $F$. The class $\W_0$ of the
\emph{K\"ahler manifolds with Norden metric} is contained in each
of the other seven classes. It is determined by the condition
$F(x,y,z)=0$, which is equivalent to $\n J=0$.

The condition                                                                        
\begin{equation}\label{1.2}
\mathop{\s} \limits_{x,y,z} F(x,y,z)=0,
\end{equation}
where $\mathop{\s} \limits_{x,y,z}$ is the cyclic sum over
$x,y,z$, characterizes the class $\W_3$ of the
\emph{quasi-K\"ahler manifolds with Norden metric}. This is the
only class of manifolds with non-integrable almost complex
structure $J$.

Let $\{e_i\}$ ($i=1,2,\dots,2n$) be an arbitrary basis of the
tangent space of $M$ at a point $p\in M$. The components of the
inverse matrix of $g$, with respect to this basis, are denoted by
$g^{ij}$.

Following \cite{GRMa}, the \emph{square norm} $\nJ$ of $\nabla J$
is defined in \cite{MeMa} by
\begin{equation}\label{1.3}
    \nJ=g^{ij}g^{ks}
    g\bigl(\left(\nabla_{e_i} J\right)e_k,\left(\nabla_{e_j}
    J\right)e_s\bigr),
\end{equation}
where it is proven that
\begin{equation}\label{1.4}
    \nJ=-2g^{ij}g^{ks}
    g\bigl(\left(\nabla_{e_i} J\right)e_k,\left(\nabla_{e_s}
    J\right)e_j\bigr).
\end{equation}

There, the manifold with $\nJ=0$ is called an
\emph{isotropic-K\"ahler manifold with Norden metric}. It is clear
that every K\"ahler manifold with Norden metric is
isotropic-K\"ahler, but the inverse implication is not always
true.

Let $R$ be the curvature tensor of $\nabla$, i.~e. $
R(x,y)z=\nabla_x \nabla_y z - \nabla_y \nabla_x z -
    \nabla_{[x,y]}z$. The corresponding $(0,4)$-tensor is
determined by $R(x,y,z,w)=g(R(x,y)z,w)$. The Ricci tensor $\rho$
and the scalar curvature $\tau$ with respect to $\nabla$ are
defined by
\[
    \rho(y,z)=g^{ij}R(e_i,y,z,e_j),\qquad \tau=g^{ij}\rho(e_i,e_j).
\]

A tensor $L$ of type (0,4) with the pro\-per\-ties%
\be{1.5}%
L(x,y,z,w)=-L(y,x,z,w)=-L(x,y,w,z), \ee %
\be{1.6} %
\mathop{\s} \limits_{x,y,z} L(x,y,z,w)=0 \quad
\textit{(the first Bianchi identity)}
\ee %
is called a \emph{curvature-like tensor}. Moreover, if the
curvature-like tensor $L$ has the property
\begin{equation}\label{1.7}
L(x,y,Jz,Jw)=-L(x,y,z,w),
\end{equation}
it is called a \emph{K\"ahler tensor} \cite{GaGrMi}.


Let $\n'$ be a linear connection with a tensor $Q$ of the
transformation $\n \rightarrow\n'$ and a torsion tensor $T$, i.~e.
\begin{equation}\label{1.8}
\n'_x y=\n_x y+Q(x,y),\qquad T(x,y)=\n'_x y-\n'_y x-[x,y].
\end{equation}
The corresponding (0,3)-tensors are defined by
\begin{equation}\label{1.9}
    Q(x,y,z)=g(Q(x,y),z), \qquad T(x,y,z)=g(T(x,y),z).
\end{equation}
The symmetry of the Levi-Civita connection implies
\begin{equation}\label{1.10}
    T(x,y)=Q(x,y)-Q(y,x), \qquad
    T(x,y)=-T(y,x).
\end{equation}

A partial decomposition of the space $\mathcal{T}$ of the torsion
(0,3)-tensors $T$ (i.~e. $T(x,y,z)=-T(y,x,z)$) is valid on an
almost complex manifold with Norden metric $(M,J,g)$:
$\mathcal{T}=\mathcal{T}_1\oplus\mathcal{T}_2\oplus\mathcal{T}_3\oplus\mathcal{T}_4$,
where $\mathcal{T}_i$ $(i=1,2,3,4)$ are invariant orthogonal
subspaces \cite{GaMi}. For the projection operators $p_i$ of
$\mathcal{T}$ in $\mathcal{T}_i$ it is established that:
\begin{equation}\label{1.11}
  \begin{array}{l}
    4p_1(x,y,z)=T(x,y,z)-T(Jx,Jy,z)-T(Jx,y,Jz)-T(x,Jy,Jz),\\[4pt]
    4p_2(x,y,z)=T(x,y,z)-T(Jx,Jy,z)+T(Jx,y,Jz)+T(x,Jy,Jz),\\[4pt]
    8p_3(x,y,z)=2T(x,y,z)-T(y,z,x)-T(z,x,y)-T(Jy,z,Jx)\\[4pt]
    \phantom{p_3(x,y,z)=-}-T(z,Jx,Jy)+2T(Jx,Jy,z)-T(Jy,Jz,x)\\[4pt]
    \phantom{p_3(x,y,z)=-}-T(Jz,Jx,y)+T(y,Jz,Jx)+T(Jz,x,Jy),\\[4pt]
    8p_4(x,y,z)=2T(x,y,z)+T(y,z,x)+T(z,x,y)+T(Jy,z,Jx)\\[4pt]
    \phantom{p_3(x,y,z)=-}+T(z,Jx,Jy)+2T(Jx,Jy,z)+T(Jy,Jz,x)\\[4pt]
    \phantom{p_3(x,y,z)=-}+T(Jz,Jx,y)-T(y,Jz,Jx)-T(Jz,x,Jy).\\[4pt]
  \end{array}
\end{equation}
A linear connection $\n'$ on an almost complex manifold with
Norden metric $(M,J,g)$ is called a \emph{natural connection} if
$\n' J=\n' g=0$. The last conditions are equivalent to $\n' g=\n'
\tilde{g}=0$. If $\n'$ is a linear connection with a tensor $Q$ of
the transformation $\n \rightarrow\n'$ on an almost complex
manifold with Norden metric, then it is  a natural connection iff
the following conditions are valid:
\begin{equation}\label{1.12}
    F(x,y,z)=Q(x,y,Jz)-Q(x,Jy,z),
\end{equation}
\begin{equation}\label{1.13}
    Q(x,y,z)=-Q(x,z,y).
\end{equation}

Since  $\n' g=0$, equalities \eqref{1.10} and \eqref{1.11} imply
\begin{equation}\label{1.14}
    Q(x,y,z)=\frac{1}{2}\bigl\{
    T(x,y,z)-T(y,z,x)+T(z,x,y)\bigr\},
\end{equation}
which is the Hayden theorem (\cite{Hay}).

\begin{thm}[\cite{Mek2}]\label{thm-1.1}
    For a natural connection with a torsion tensor $T$ on a quasi-K\"ahler manifold with Norden metric
    $(M,J,g)$, which is non-K\"ahlerian, the following properties
    are valid
    \begin{equation}\label{3.12}
        p_2\neq 0,\qquad p_3=0.
    \end{equation}\hfill$\Box$
\end{thm}


\section{Connection with totally skew-symmetric torsion}

\subsection{Main results}

The main results in the present work are consist in the following
\begin{thm}\label{thm-2.1}
Let $\n'$ be a natural connection with a totally skew-symmetric
torsion tensor on a non-K\"ahler almost complex manifold with
Norden metric $(M,J,g)$. Then
\begin{enumerate}    \renewcommand{\labelenumi}{(\roman{enumi})}
    \item $(M,J,g)$ is quasi-K\"ahlerian;
    \item $\n'$ has a torsion tensor $T$ in the class
    $\mathcal{T}_2\oplus\mathcal{T}_4$;
    \item $\n'$ is unique and the
    tensor $Q$ of the transformation $\n \rightarrow\n'$ is
    determined by
    \begin{equation}\label{2.1}
        Q(x,y,z)=\frac{1}{4}\left\{F(x,Jy,z)-F(Jx,y,z)-2F(y,Jx,z)\right\}.
    \end{equation}
\end{enumerate}
\end{thm}

\begin{proof}
Since $\n'$ has a totally skew-symmetric torsion tensor $T$, then
we have
\begin{equation}\label{2.2}
    T(x,y,z)=-T(y,x,z)=-T(x,z,y)=-T(z,y,x).
\end{equation}
From \eqref{1.14} and \eqref{2.2} it is follows that for the
tensor $Q$ of the transformation $\n \rightarrow\n'$ it is valid
\begin{equation}\label{2.3}
    Q(x,y,z)=\frac{1}{2}T(x,y,z).
\end{equation}
Since $\n'$ is a natural connection, then \eqref{1.12} holds and
consequently
\begin{equation}\label{2.4}
    \mathop{\s} \limits_{x,y,z} F(x,y,z)=\mathop{\s} \limits_{x,y,z}\left\{Q(x,y,Jz)-Q(x,Jy,z)\right\}.
\end{equation}
According to \eqref{2.3}, $Q$ is also 3-form and then we have
\begin{equation}\label{2.5}
    Q(x,y,Jz)=Q(y,Jz,x).
\end{equation}

Equalities \eqref{2.4} and \eqref{2.5} imply \eqref{1.2}, which
completes the proposition (i).

According to \thmref{thm-1.1}, for the tensor $T$ we have $p_2\neq
0$ and $p_3=0$, i.~e.
$T\in\mathcal{T}_1\oplus\mathcal{T}_2\oplus\mathcal{T}_4$. Let us
suppose that $T\in\mathcal{T}_1$. Then $T=p_1$ and from
\eqref{1.11} it is follows
\begin{equation}\label{2.6}
    3T(x,y,z)+T(Jx,Jy,z)+T(Jx,y,Jz)+T(x,Jy,Jz)=0.
\end{equation}
We substitute $Jy$ for $y$ and $Jz$ for $z$ in \eqref{2.6} and the
obtained equality we add to \eqref{2.6}. In the result we
substitute $Jz$ for $z$ and according to \eqref{2.3} we get
$Q(x,y,Jz)-Q(x,Jy,z)=0$. Then, by \eqref{1.12} we obtain
$F(x,y,z)=0$, i.~e. $(M,J,g)$ is K\"ahlerian, which is a
contradiction. Therefore $T\notin\mathcal{T}_1$, i.~e. $p_1=0$.

Combining \eqref{1.11}, \eqref{2.2} and \eqref{2.3}, we find that
\begin{equation}\label{2.7}
    p_4(x,y,z)=Q(x,y,z)+Q(Jx,Jy,z).
\end{equation}
Let us suppose that $p_4=0$. Then from \eqref{2.7} we have
\begin{equation}\label{2.8}
    Q(x,y,z)+Q(Jx,Jy,z)=0.
\end{equation}
We substitute $y  \leftrightarrow z$ in \eqref{2.8}, and subtract
the obtained equality from \eqref{2.8}. In the result we apply
\eqref{2.8} and the fact that $Q$ is a 3-form. In such a way we
obtain
\begin{equation}\label{2.9}
    Q(x,y,z)+Q(Jx,y,Jz)=0.
\end{equation}
The equalities \eqref{2.8} and \eqref{2.9} imply
$Q(x,y,Jz)-Q(x,Jy,z)=0$. Thus, according to \eqref{1.12}, we
obtain $F(x,y,z)=0$, i.~e. $(M,J,g)$ is K\"ahlerian, which is
impossible. Therefore $p_4\neq 0$.

Thereby we establish the following conditions for $T$: $p_1=0$,
$p_2\neq 0$, $p_3=0$, $p_4\neq 0$. Hence
$T\in\mathcal{T}_2\oplus\mathcal{T}_4$, i.~e. (ii) holds.

From \eqref{1.12}, having in mind that $Q$ is a 3-form, we obtain
\begin{equation}\label{2.10}
\begin{array}{l}
    F(x,Jy,z)-F(Jx,y,z)-2F(y,Jx,z)\\[4pt]
    =Q(Jx,Jy,z)+Q(Jx,y,Jz)+Q(x,Jy,Jz)+3Q(x,y,z).
\end{array}
\end{equation}
Since $p_1=0$, from \eqref{1.11} and \eqref{2.3} we have
\begin{equation}\label{2.11}
    Q(Jx,Jy,z)+Q(Jx,y,Jz)+Q(x,Jy,Jz)=Q(x,y,z).
\end{equation}
The equalities \eqref{2.10} and \eqref{2.11} imply \eqref{2.1},
which completes the proof of (iii).
\end{proof}

According to \eqref{1.1}, equality \eqref{2.1} is equivalent to
\begin{equation}\label{2.12}
    Q(x,y)=\frac{1}{4}\bigl\{\left(\n_{x} J\right)Jy-\left(\n_{Jx} J\right)y-2\left(\n_{y} J\right)Jx
    \bigr\}.
\end{equation}

\begin{rem}
In \cite{Teo} it is found a 2-parametric family of natural
connections on almost complex manifolds with Norden metric which
contains the connection $\n'$ with totally skew-symmetric torsion
tensor.
\end{rem}

\subsection{An example}
\label{sec-exa}

Let $V$ be a real 4-dimensional vector space with a basis
$\{E_i\}$. Let us consider a structure of a Lie algebra determined
by the commutators $[E_i,E_j]=C_{ij}^k E_k$, where $C_{ij}^k$ are
structure constants satisfying the anti-commutativity condition
$C_{ij}^k=-C_{ji}^k$ and the Jacobi identity $C_{ij}^k
C_{ks}^l+C_{js}^k C_{ki}^l+C_{si}^k C_{kj}^l=0$.

Let $G$ be the associated connected Lie group and $\{X_i\}$ be a
global basis for the left invariant vector fields that is induced
by the basis $\{E_i\}$ of $V$. Then we have the decomposition
\begin{equation}\label{2.0}
    [X_i,X_j]=C_{ij}^k X_k.
\end{equation}

Let us consider the almost complex manifold with Norden metric
$(G,J,g)$, where
\begin{equation}\label{2.13}
    JX_1=X_3,\qquad JX_2=X_4,\qquad JX_3=-X_1,\qquad JX_4=-X_2
\end{equation}
and
\begin{equation}\label{2.14}
\begin{array}{c}
  g(X_1,X_1)=g(X_2,X_2)=-g(X_3,X_3)=-g(X_4,X_4)=1, \\[4pt]
  g(X_i,X_j)=0\quad \text{for}\quad i\neq j. \\
\end{array}
\end{equation}

Because of \eqref{2.14} the following equality is valid
\begin{equation}\label{2.15}
    2g\left(\n_{X_i} X_j,X_k
    \right)=g\left([X_i,X_j],X_k\right)+g\left([X_k,X_i],X_j\right)+g\left([X_k,X_j],X_i\right).
\end{equation}

We add the condition for the associated metric $\tilde{g}$ of $g$
to be a Killing metric \cite{Hel}. Then
\begin{equation}\label{2.16}
    g\left([X_i,X_j],JX_k\right)+g\left([X_i,X_k],JX_j\right)=0.
\end{equation}

Combining \eqref{1.1}, \eqref{2.13}, \eqref{2.14}, \eqref{2.15}
and \eqref{2.16}, we obtain
\begin{equation}\label{2.17}
    2F(X_i,X_j,X_k)=g\left([JX_i,X_j],X_k\right)+g\left([JX_i,X_k],X_j\right).
\end{equation}
From \eqref{2.17} it is follows immediately that $\mathop{\s}
\limits_{X_i,X_j,X_k} F(X_i,X_j,X_k)=0$, i.~e. $(G,J,g)$ is a
quasi-K\"ahler manifold with Norden metric.

According to \eqref{2.0} and \eqref{2.16}, we get
\begin{equation}\label{2.18}
\begin{array}{ll}
    [X_1,X_2]= \lambda_1 X_1 +\lambda_2 X_2,\qquad & [X_1,X_3]= \lambda_3
    X_2-\lambda_1 X_4,\\[4pt]
    [X_1,X_4]= -\lambda_3 X_1 -\lambda_2 X_4,\qquad & [X_2,X_3]= \lambda_4
    X_2+\lambda_1 X_3,\\[4pt]
    [X_2,X_4]= -\lambda_4 X_1 +\lambda_2 X_3,\qquad & [X_3,X_4]= \lambda_3
    X_3
    +\lambda_4 X_4,\\[4pt]
\end{array}
\end{equation}
where
\[
\begin{array}{ll}
    \lambda_1=C_{12}^1=C_{23}^3=-C_{13}^4, \qquad & \lambda_2=C_{12}^2=C_{24}^3=-C_{14}^4,\\[4pt]
    \lambda_3=C_{13}^2=C_{34}^3=-C_{14}^1, \qquad & \lambda_4=C_{23}^2=C_{34}^4=-C_{24}^1.\\[4pt]
\end{array}
\]

Let equalities \eqref{2.18} are valid for an almost complex
manifold with Norden metric $(G,J,g)$, where $J$ and $g$ are
determined by \eqref{2.13} and \eqref{2.14}. Then we verify
directly that  the Jacobi identity for the commutators $[X_i,X_j]$
is satisfied and the associated metric $\tilde{g}$ of $g$ is a
Killing metric.

Therefore, the following theorem is valid.
\begin{thm}\label{thm-2.2}
    Let $(G,J,g)$ be a 4-dimensional
almost complex manifold with Norden metric, where $G$ is the
connected Lie group with an associated Lie algebra determined by a
global basis $\{X_i\}$ of left invariant vector fields, and $J$
and $g$ are the almost complex structure and the Norden metric
determined by \eqref{2.13} and \eqref{2.14}, respectively. Then
$(G,J,g)$ is a quasi-K\"ahler manifold with a Killing associated
metric $\tilde{g}$ iff $G$ belongs to the 4-parametric family of
Lie groups, defined by \eqref{2.18}.\hfill$\Box$
\end{thm}

Let $(G,J,g)$ be the quasi-K\"ahler manifold determined by the
conditions of \thmref{thm-2.2}.

By \eqref{2.17} and \eqref{2.18} we get the non-trivial components
$F_{ijk}=F(X_i,X_j,X_k)$ of the tensor $F$:
\[
\begin{array}{ll}
    -2F_{114}=2F_{123}=-2F_{312}=2F_{334}=2F_{411}=2F_{433}=\lambda_1, \\[4pt]
    2F_{223}=-2F_{241}=-2F_{322}=-2F_{344}=2F_{412}=2F_{434}=\lambda_2,\\[4pt]
    -2F_{112}=-2F_{134}=F_{211}=F_{233}=2F_{314}=-2F_{332}=\lambda_3, \\[4pt]
    -F_{122}=F_{144}=2F_{221}=2F_{234}=2F_{414}=-2F_{432}=\lambda_4.\\[4pt]
\end{array}
\]

This leads to the following
\begin{prop}\label{prop-2.3}
$(G,J,g)$ is a non-K\"ahler manifold with Norden
metric.\hfill$\Box$
\end{prop}

Using \eqref{2.15} and \eqref{2.16}, we obtain
\[
    2g\left(\n_{X_i} X_j,X_k
    \right)=g\left([X_i,X_j],X_k\right)-g\left(J[X_i,JX_j],X_k\right)+g\left(J[X_i,X_j],X_k\right).
\]
Then we have
\begin{equation}\label{2.19}
   2 \n_{X_i} X_j
    =[X_i,X_j]-J[X_i,JX_j]+J[JX_i,X_j].
\end{equation}

The equality \eqref{2.19} implies
\[
\begin{array}{l}
    2\n_{X_i} JX_j=[X_i,JX_j]+J[X_i,X_j]+J[JX_i,JX_j],\\[4pt]
    2J\n_{X_i} X_j=J[X_i,X_j]+[X_i,JX_j]-[JX_i,X_j].\\[4pt]
\end{array}
\]
We subtract the last two equalities, apply the formula for
covariant derivation and obtain
\begin{equation}\label{2.20}
       2 \left(\n_{X_i} J\right)X_j=J[JX_i,JX_j]+[JX_i,X_j].
\end{equation}

By \eqref{2.18} and \eqref{2.20} we obtain the components of $\n
J$:
\begin{equation}\label{2.21}
\begin{array}{l}
        2\left(\n_{X_1} J\right)X_1=2\left(\n_{X_3} J\right)X_3=-\lambda_3 X_2 + \lambda_1 X_4,\\[4pt]
        2\left(\n_{X_2} J\right)X_2=2\left(\n_{X_4} J\right)X_4= \lambda_4 X_1 - \lambda_2 X_3,\\[4pt]
        2\left(\n_{X_1} J\right)X_3=-2\left(\n_{X_3} J\right)X_1= \lambda_1 X_2 + \lambda_3 X_4,\\[4pt]
        2\left(\n_{X_2} J\right)X_4=-2\left(\n_{X_4} J\right)X_2=-\lambda_4 X_3 - \lambda_2 X_1,\\[4pt]
        2\left(\n_{X_1} J\right)X_2=-\lambda_3 X_1 -2\lambda_4 X_2 -\lambda_1 X_3,\\[4pt]
        2\left(\n_{X_1} J\right)X_4=-\lambda_1 X_1 +\lambda_3 X_3 +2\lambda_4 X_4,\\[4pt]
        2\left(\n_{X_2} J\right)X_1=2\lambda_3 X_1 +\lambda_4 X_2 +\lambda_2 X_4,\\[4pt]
        2\left(\n_{X_2} J\right)X_3= \lambda_2 X_2 -2\lambda_3 X_3 -\lambda_4 X_4,\\[4pt]
        2\left(\n_{X_3} J\right)X_2=-\lambda_1 X_1 -2\lambda_2 X_2 +\lambda_3 X_3,\\[4pt]
        2\left(\n_{X_3} J\right)X_4= \lambda_3 X_1 +\lambda_1 X_3 +2\lambda_2 X_4,\\[4pt]
        2\left(\n_{X_4} J\right)X_1=2\lambda_1 X_1 +\lambda_2 X_2 -\lambda_4 X_4,\\[4pt]
        2\left(\n_{X_4} J\right)X_3=-\lambda_4 X_2 -2\lambda_1 X_3 -\lambda_2 X_4.\\[4pt]
\end{array}
\end{equation}

From \eqref{1.3} and  \eqref{2.21} it follows that
$\nJ=-4\left(\lambda_1^2+\lambda_2^2-\lambda_3^2-\lambda_4^2\right)$,
and therefore the following is valid.
\begin{prop}\label{prop-2.4}
$(G,J,g)$ is an isotropic K\"ahler manifold with Norden metric iff
$\lambda_1^2+\lambda_2^2-\lambda_3^2-\lambda_4^2=0$. \hfill$\Box$
\end{prop}

Let $\n'$ be the connection with totally skew-symmetric torsion
tensor $T$ on $(G,J,g)$.

By virtue of \eqref{2.3}, \eqref{2.12} and  \eqref{2.20}, we have
\begin{equation}\label{2.22}
        4T(X_i,X_j)=3[JX_i,JX_j]-J[X_i,JX_j]-J[JX_i,X_j]+[X_i,X_j],
\end{equation}
which implies
\begin{equation}\label{2.23}
    \begin{array}{ll}
    T(X_1,X_2)= \lambda_3 X_3 +\lambda_4 X_4,\qquad & T(X_1,X_3)=
    \lambda_3 X_2-\lambda_1 X_4,\\[4pt]
    T(X_1,X_4)= \lambda_4 X_2 +\lambda_1 X_3,\qquad & T(X_2,X_3)=
    -\lambda_3 X_1-\lambda_2 X_4,\\[4pt]
    T(X_2,X_4)= -\lambda_4 X_1 +\lambda_2 X_3,\qquad & T(X_3,X_4)=
    \lambda_1 X_1 +\lambda_2 X_2.\\[4pt]
\end{array}
\end{equation}

Then the non-trivial components $T_{ijk}=T(X_i,X_j,X_k)$ of the
corresponding (0,3)-tensor $T$ are:
\begin{equation}\label{2.24}
    T_{134}=\lambda_1,\quad T_{234}=\lambda_2,\quad T_{123}=-\lambda_3,\quad T_{124}=-\lambda_4.
\end{equation}

According to \eqref{2.18} and  \eqref{2.23}, we have
\[
T(X_i,X_j)=[JX_i,JX_j]
\]
and then, from \eqref{2.22} it follows that
\begin{equation}\label{2.25}
    [JX_i,JX_j]+J[JX_i,X_j]=[X_i,X_j]-J[X_i,JX_j].
\end{equation}

Since we have $2Q=T$ for $\n'$, from \eqref{1.8}, \eqref{2.22} and
\eqref{2.25}, we obtain
\begin{equation}\label{2.26}
    \n'_{X_i} X_j=[X_i,X_j]-J[X_i,JX_j].
\end{equation}

Combining \eqref{2.18} and  \eqref{2.26}, we get the components of
$\n'$:
\begin{equation}\label{2.27}
    \begin{array}{l}
        \n'_{X_1} X_1=\n'_{X_3} X_3= -\lambda_1 X_2 -\lambda_3 X_4,\\[4pt]
        \n'_{X_2} X_2=\n'_{X_4} X_4=  \lambda_2 X_1 +\lambda_4 X_3,\\[4pt]
        \n'_{X_1} X_2=\n'_{X_3} X_4=  \lambda_1 X_1 +\lambda_3 X_3,\\[4pt]
        \n'_{X_1} X_3=-\n'_{X_3} X_1= \lambda_3 X_2 -\lambda_1 X_4,\\[4pt]
        \n'_{X_1} X_4=-\n'_{X_3} X_2= -\lambda_3 X_1 +\lambda_1 X_3,\\[4pt]
        \n'_{X_2} X_1=\n'_{X_4} X_3= -\lambda_2 X_2 -\lambda_4 X_4,\\[4pt]
        \n'_{X_2} X_3=-\n'_{X_4} X_1= \lambda_4 X_2 -\lambda_2 X_4,\\[4pt]
        \n'_{X_2} X_4=-\n'_{X_4} X_2= -\lambda_4 X_1 +\lambda_2 X_3.\\[4pt]
\end{array}
\end{equation}

Thus we arrive at the following
\begin{prop}\label{prop-2.5}
Let $\n'$ be the connection with totally skew-symmetric torsion
tensor $T$ on $(G,J,g)$. Then the components of $\n'$ and $T$ with
respect to the basis $\{X_i\}$ are \eqref{2.27} and \eqref{2.24},
respectively. \hfill$\Box$
\end{prop}

In \cite{Mek2} and \cite{Mek3} we have considered the canonical
connection $\n^C$ and the B-connection $\n^B$  on quasi-K\"ahler
manifolds with Norden metric. Now we will show how these two
connections relate to the connection $\n'$ with totally
skew-symmetric torsion tensor on $(G,J,g)$.

The connection $\n^B$ is defined in \cite{Mek3} by
\[
\n^B_{X_i} X_j=\n_{X_i} X_j+\frac{1}{2}\left(\n_{X_i}J\right)JX_j.
\]
Hence, applying \eqref{2.13} and \eqref{2.21}, we obtain directly
$\n^B=\frac{3}{4}\n'$. On the other hand,  according to
\cite{Mek2}, we have $\n^B=\frac{1}{2}\left(\n^C+\n'\right)$. Thus
we proved the following
\begin{prop}\label{prop-2.6}
Let $\n'$, $\n^B$ and $\n^C$ be the connection with totally
skew-symmetric torsion, the B-connection and the canonical
connection on $(G,J,g)$, respectively. Then
\[
\n'=\frac{4}{3}\n^B=2\n^C.
\]\hfill$\Box$
\end{prop}

In \cite{Mek1} we have proved that $\n'$ has a K\"ahler curvature
tensor on any quasi-K\"ahler manifolds with Norden metric  iff the
following identity holds
\begin{equation}\label{2.28}
       \mathop{\s} \limits_{x,y,z} \bigl\{
      g\bigl(\left(\nabla_x J\right)Jy+\left(\nabla_{Jx} J\right)y,
      \left(\nabla_z
       J\right)Jw+\left(\nabla_{Jz}J\right)w\bigr)
    \bigr\}=0.
\end{equation}

Equality \eqref{2.20} implies
\[
\begin{array}{l}
    2\left(\n_{X_i} J\right)JX_j=-J[JX_i,X_j]+[JX_i,JX_j],\\[4pt]
    2\left(\n_{JX_i} J\right)X_j=-J[X_i,JX_j]-[X_i,X_j].\\[4pt]
\end{array}
\]
We add the above, and according to \eqref{2.25} obtain
\[
    \left(\nabla_{X_i} J\right)JX_j+\left(\nabla_{JX_i} J\right)X_j
    =
      [JX_i,JX_j]-[X_i,X_j].
\]
Then, the condition \eqref{2.28} takes the following form on
$(G,J,g)$:
\[
       \mathop{\s} \limits_{x,y,z} \bigl\{
       g\bigl([JX_i,JX_j]-[X_i,X_j],
      [JX_k,JX_s]-[X_k,X_s]\bigr)\bigr\}=0.
\]
The last equality and equalities \eqref{2.13},  \eqref{2.14} and
\eqref{2.18} imply
\begin{prop}\label{prop-2.7}
The connection $\n'$ with totally skew-symmetric torsion  on
$(G,J,g)$ has a K\"ahler curvature tensor iff
$\lambda_1^2+\lambda_2^2=\lambda_3^2+\lambda_4^2$. \hfill$\Box$
\end{prop}


\section{Connection with parallel totally skew-sym\-met\-ric torsion}

Let $\n'$ be the connection with totally skew-symmetric torsion
tensor $T$ on the quasi-K\"ahler manifold with Norden metric
$(M,J,g)$.

Now we consider the case  when $\n'$ has a parallel torsion, i.~e.
$\n'T=0$.

It is known that the curvature tensors $R'$ and $R$ of $\n'$ and
$\n$, respectively,  satisfy:
\begin{equation}\label{3.1}
    \begin{split}
    &R'(x,y,z,w)=R(x,y,z,w)+ \left(\n_x Q\right)(y,z,w)-\left(\n_y
    Q\right)(x,z,w)
    \\[4pt]
    &\phantom{K(x,y,z,w)=R(x,y,z,w)}+Q\bigl(x,Q(y,z),w\bigr)-Q\bigl(y,Q(x,z),w\bigr).
\end{split}
\end{equation}

Equality \eqref{2.3} implies $\n'Q=0$ in the considered case. Then
from the formula for covariant derivation with respect to $\n'$ it
follows that
\begin{equation}\label{3.2}
    xQ(y,z,w)-Q(\n'_x y,z,w)-Q(y,\n'_x z,w)-Q(y,z,\n'_x w)=0.
\end{equation}
According to the first equality of \eqref{1.8} we have
\begin{equation}\label{3.3}
    \begin{array}{l}
Q(\n'_x y,z,w)=Q(\n_x y,z,w)+Q(Q(x, y),z,w),\\[4pt]
Q( y,\n'_x z,w)=Q(y,\n_x z,w)+Q(y, Q(x, z),w),\\[4pt]
Q( y,z,\n'_x w)=Q(y,z,\n_x w)+Q(y,  z,Q(x,w)).\\[4pt]
    \end{array}
\end{equation}
Combining \eqref{3.2}, \eqref{3.3}, the first equality of
\eqref{1.9} and having in mind the formula for covariant
derivation with respect to $\n$, we obtain
\begin{equation}\label{3.4}
    \begin{split}
    &\left(\n_x Q\right)(y,z,w)=Q(Q(x, y),z,w)    \\[4pt]
    &\phantom{\left(\n_x Q\right)(y,z,w)=}
    -g\bigl(Q(x,z),Q(y,w)\bigr)-g\bigl(Q(y,z),Q(x,w)\bigr).
\end{split}
\end{equation}
From \eqref{3.4} and the first equality of \eqref{1.10} we have
\begin{equation}\label{3.5}
    \begin{split}
    &\left(\n_x Q\right)(y,z,w)-\left(\n_y Q\right)(x,z,w)
    =Q(T(x, y),z,w)    \\[4pt]
    &\phantom{\left(\n_x Q\right)(y,z,w)=}
    -2g\bigl(Q(x,z),Q(y,w)\bigr)+2g\bigl(Q(y,z),Q(x,w)\bigr).
\end{split}
\end{equation}

Because of \eqref{3.5}, equality \eqref{3.1} can be rewritten as
\begin{equation}\label{3.6}
    \begin{split}
    &R'(x,y,z,w)=R(x,y,z,w)+ Q\left(T(x,y),z,w\right)
    \\[4pt]
    &\phantom{R'(x,y,z,w)=}
    -g\bigl(Q(x,z),Q(y,w)\bigr)+g\bigl(Q(y,z),Q(x,w)\bigr).
    \end{split}
\end{equation}

Since $Q(e_i,e_j)=-Q(e_j,e_i)$ it follows that
$g^{ij}Q(e_i,e_j)=0$. Then, from \eqref{3.6} after contraction by
$x=e_i$, $w=e_j$, we obtain the following equality for the Ricci
tensor $\rho'$ of $\n'$:
\begin{equation}\label{3.7}
    \begin{split}
    &\rho'(y,z)=\rho(y,z)+2g^{ij}g\bigl(Q(e_i,y),Q(z,e_j)\bigr)
    \\[4pt]
    &\phantom{\rho'(y,z)=\rho(y,z)}
    -g^{ij}g\bigl(Q(e_i,z),Q(y,e_j)\bigr).
    \end{split}
\end{equation}
Contracting by $y=e_k$, $z=e_s$ in \eqref{3.7}, we get
\begin{equation}\label{3.8}
    \tau'=\tau+g^{ij}g^{ks}g\bigl(Q(e_i,e_k),Q(e_s,e_j)\bigr),
\end{equation}
where $\tau'$ is the scalar curvature of $\n'$.

By virtue of \eqref{3.8}, \eqref{2.12}, \eqref{1.3} and
\eqref{1.4} we have
\begin{equation}\label{3.9}
    \tau'=\tau-\frac{1}{8}\nJ.
\end{equation}

Thus we arrive at the following
\begin{thm}\label{thm-3.1}
Let $\n'$ be the connection with parallel totally skew-symmetric
torsion on the quasi-K\"ahler manifold with Norden metric
$(M,J,g)$. Then for the Ricci tensor $\rho'$ and the scalar
curvature $\tau'$ of $\n'$ are valid \eqref{3.7} and \eqref{3.9},
respectively. \hfill$\Box$
\end{thm}

Equality \eqref{3.9} leads to the following
\begin{cor}\label{cor-3.2}
Let $\n'$ be the connection with parallel totally skew-symmetric
torsion on the quasi-K\"ahler manifold with Norden metric
$(M,J,g)$. Then the manifold $(M,J,g)$ is isotropic K\"ahlerian
iff $\n'$ and $\n$ have equal scalar curvatures. \hfill$\Box$
\end{cor}

\begin{rem}
The  4-parametric family of 4-dimensional quasi-K\"ahler manifolds
$(G,J,g)$ considered in subsection \ref{sec-exa} does not admit
any connection with parallel totally skew-symmetric torsion.
\end{rem}


\section{Connection with parallel totally skew-sym\-met\-ric torsion and K\"ahler curvature tensor}

Let $\n'$ be a connection with parallel totally skew-symmetric
torsion on the quasi-K\"ahler manifold with Norden metric
$(M,J,g)$.

We will find conditions for the curvature tensor $R'$ of $\n'$ to
be K\"ahlerian.

From \eqref{2.3}, having in mind that $Q$ is a 3-form, we have
\begin{equation*}
    \begin{split}
    &Q\left(T(x,y),z,w\right)=Q\left(z,w,T(x,y)\right)=g\bigl(Q(z,w),T(x,y)\bigr)
    \\[4pt]
    &\phantom{Q\left(T(x,y),z,w\right)}
    =g\bigl(T(x,y),Q(z,w)\bigr)=2g\bigl(Q(x,y),Q(z,w)\bigr).
    \end{split}
\end{equation*}

Then \eqref{3.6} obtains the form
\begin{equation}\label{4.1}
    \begin{split}
    &R'(x,y,z,w)=R(x,y,z,w)+2g\bigl(Q(x,y),Q(z,w)\bigr)
    \\[4pt]
    &\phantom{R'(x,y,z,w)=}
    -g\bigl(Q(x,z),Q(y,w)\bigr)+g\bigl(Q(y,z),Q(x,w)\bigr).
    \end{split}
\end{equation}

From \eqref{4.1}, identities \eqref{1.5} and \eqref{1.7} for $R'$
follow immediately. Therefore $R'$ is a K\"ahler tensor iff the
first Bianchi identity \eqref{1.6} for $R'$ is satisfied. Since
this identity is valid for $R$, then \eqref{4.1} implies that $R'$
is K\"ahlerian iff
\[
\mathop{\s} \limits_{x,y,z}
\bigl\{2g\bigl(Q(x,y),Q(z,w)\bigr)-g\bigl(Q(x,z),Q(y,w)\bigr)+g\bigl(Q(y,z),Q(x,w)\bigr)\bigr\}=0.
\]

Thus, using that $Q$ is a skew-symmetric tensor, we arrive the
following
\begin{thm}\label{thm-4.1}
Let $\n'$ be the connection with parallel totally skew-symmetric
torsion on the quasi-K\"ahler manifold with Norden metric
$(M,J,g)$. Then the curvature tensor for $\n'$ is a K\"ahler
tensor iff
\begin{equation}\label{4.2}
    \mathop{\s} \limits_{x,y,z}
g\bigl(Q(x,y),Q(z,w)\bigr)=0.
\end{equation} \hfill$\Box$
\end{thm}

Because of the skew-symmetry of $Q$, \eqref{4.2} implies
\begin{equation*}
    g\bigl(Q(y,z),Q(x,w)\bigr)-g\bigl(Q(x,z),Q(y,w)\bigr)=-g\bigl(Q(x,y),Q(z,w)\bigr).
\end{equation*}
The last equality and \eqref{4.1} lead to the following
\begin{cor}\label{cor-4.2}
Let $\n'$ be the connection with parallel totally skew-symmetric
torsion and K\"ahler curvature tensor on the quasi-K\"ahler
manifold with Norden metric $(M,J,g)$. Then
\begin{equation}\label{4.3}
    R'(x,y,z,w)=R(x,y,z,w)+g\bigl(Q(x,y),Q(z,w)\bigr).
\end{equation}
\hfill$\Box$
\end{cor}

If $R'$ is a K\"ahler tensor then $R'(x,y,Jz,Jw)=-R'(x,y,z,w)$,
and because of \eqref{4.3} we have
\begin{equation}\label{4.4}
    \begin{split}
    &R(x,y,Jz,Jw)+R(x,y,z,w)=-g\bigl(Q(x,y),Q(Jz,Jw)\bigr)\\[4pt]
   &\phantom{R(x,y,Jz,Jw)+R(x,y,z,w)=}
    -g\bigl(Q(x,y),Q(z,w)\bigr).
    \end{split}
\end{equation}

From \eqref{2.12} we get
\[
Q(x,Jy)=JQ(x,y)-\left(\n_x J\right)y.
\]
Then we have
\[
Q(Jx,Jy)=-Q(x,y)-\left(\n_{Jx} J\right)y-\left(\n_{y} J\right)Jx
\]
and consequently
\[
    \begin{split}
&g\bigl(Q(x,y),Q(Jz,Jw)\bigr)=-g\bigl(Q(x,y),Q(z,w)\bigr)\\[4pt]
   &\phantom{g\bigl(Q(x,y),Q(Jz,Jw)\bigr)=}
-g\bigl(Q(x,y),\left(\n_{Jz} J\right)w+\left(\n_{w}
J\right)Jz\bigr).
    \end{split}
\]

The last equality and \eqref{4.4} imply the following
\begin{cor}\label{cor-4.3}
Let $\n'$ be the connection with parallel totally skew-symmetric
torsion and K\"ahler curvature tensor on the quasi-K\"ahler
manifold with Norden metric $(M,J,g)$. Then
\begin{equation}\label{4.5}
    R(x,y,Jz,Jw)+R(x,y,z,w)=g\bigl(Q(x,y),\left(\n_{Jz} J\right)w+\left(\n_{w}
J\right)Jz\bigr).
\end{equation}
\hfill$\Box$
\end{cor}

Contracting by $x=e_i$, $w=e_j$ in \eqref{4.5}, we obtain
\[
    g^{ij}R(e_i,y,Jz,Je_j)+\rho(y,z)=g^{ij}g\bigl(Q(e_i,y),\left(\n_{Jz} J\right)e_j+\left(\n_{e_j}
J\right)Jz\bigr).
\]
Then, after a contraction by $y=e_k$, $z=e_s$, it follows
\begin{equation}\label{4.6}
        \tau^{**}+\tau=g^{ij}g^{ks}g\bigl(Q(e_i,e_k),\left(\n_{Je_s} J\right)e_j+\left(\n_{e_j}
J\right)Je_s\bigr),
\end{equation}
where $\tau^{**}=g^{ij}g^{ks}R(e_i,e_k,Je_s,Je_j)$.

From \eqref{2.12}, \eqref{1.3} and  \eqref{1.4} we have
\[
g^{ij}g^{ks}g\bigl(Q(e_i,e_k),\left(\n_{Je_s}
J\right)e_j+\left(\n_{e_j} J\right)Je_s\bigr)=-\frac{1}{8}\nJ.
\]
Then \eqref{4.6} can be rewritten as
\[
\tau^{**}+\tau=-\frac{1}{8}\nJ.
\]
On the other hand, according to  \cite{MeMa}, we have
\[
\tau^{**}+\tau=-\frac{1}{2}\nJ.
\]
Then $\nJ=0$ and therefore the following is valid.
\begin{thm}
Let $(M,J,g)$ be  a quasi-K\"ahler manifold with Norden metric
which admit a connection with parallel totally skew-symmetric
torsion and K\"ahler curvature tensor. Then $(M,J,g)$ is a
isotropic K\"ahler manifold with Norden metric. \hfill$\Box$
\end{thm}



\begin{thebibliography}{99}



\bibitem{Stro}
A. Strominger, Superstrings with torsion. \emph{Nuclear Phys. B},
274:253-284 (1986).

\bibitem{Bis}
J.-M. Bismut, A local index theorem for a non-Kahler manifolds.
\emph{Math. Ann.} 284: 681-699 (1989).

\bibitem{Gaud}
P. Gauduchon, Hermitian connections and Dirac operators.
\emph{Boll. Unione Mat. Ital. Sez. B Artic. Ric. Mat. (8)} 11:
257-289 (1997).

\bibitem{Mek1}
D. Mekerov, A connection with skew symmetric torsion and Kahler
curvature tensor on a quasi-Kahler manifolds with Norden meric.
\emph{C. R. Acad. Bulgare Sci.} 61: 1249-1256 (2008).

\bibitem{GrMeDj}
K. I. Gribachev, D. G. Mekerov and G. D. Djelepov, Generalized
B-mani\-fold. \emph{C. R. Acad. Bulgare Sci.} 38: 299-302 (1985).

\bibitem{GaBo}
G. T. Ganchev and A. V. Borisov, Note on the almost complex
manifolds with a Norden metric. \emph{C. R. Acad. Bulgare Sci.}
39: 31-34 (1986).

\bibitem{GRMa}
E. Garcia-Rio and Y. Matsushita, Isotropic Kahler structures on
Engel 4-mani\-folds. \emph{J. Geom. Phys.} 33: 288-294 (2000).

\bibitem{MeMa}
D. Mekerov and M. Manev, On the geometry of quasi-Kahler manifolds
with Norden meric. \emph{Nihonkai Math. J.} 16: 89-93 (2005).

\bibitem{GaGrMi}
G. Ganchev, K. Gribachev and V. Mihova, B-connection and their
conformal invariants on conformally Kahler manifolds with
B-metric. \emph{Publ. Inst. Math. (Beograd) (N.S.)} 42: 107-121
(1987).

\bibitem{GaMi}
G. Ganchev and V. Mihova, Canonical connection and the canonical
conformal group on an almost complex manifold with B-metric.
\emph{Annuaire Univ. Sofia Fac. Math. Inform.} 81: 195-206 (1987).

\bibitem{Hay}
H. Hayden, Subspaces of a space whit torsion. \emph{Proc. London
Math. Soc. (3)} 34: 27-50 (1934).

\bibitem{Mek2}
D. Mekerov. Canonical connection on quasi-Kahler manifolds with
Norden meric (to appear). \emph{arXiv:0812.3516} (December 2008).

\bibitem{Teo}
M. Teofilova, On the geometry of almost complex manifolds with a
Norden metric. Ph. D. Thesis, Univ. of  Plovdiv, Sofia, 2008.

\bibitem{Hel}
S. Helgason, \emph{Differential geometry, Lie groups and symmetric
spaces}. Academic Press, NY, 1978.

\bibitem{Mek3}
D. Mekerov, On the geometry of the B-connection on quasi-Kahler
manifolds with Norden meric. \emph{C. R. Acad. Bulgare Sci.} 61:
1105-1110 (2008).




\end{thebibliography}
\end{document}